\documentclass[12pt,leqno]{amsart}
\usepackage{amssymb}
\usepackage{amsxtra}

\theoremstyle{plain}
\newtheorem{Thm}{Theorem}
\newtheorem{Cor}{Corollary}

\theoremstyle{definition}
\newtheorem{Def}{Definition}
\theoremstyle{remark}

\newtheorem{Ex}{Example}

\numberwithin{equation}{section}

\newcommand{\R}{\mathbf{R}}
\newcommand{\Z}{\mathbf{Z}}
\newcommand{\sinc}{\mathrm{sinc}}
\newcommand{\sign}{\mathrm{sign}}

\begin{document}

\title[Sum of Non-Identical Uniform Random Variables]{On the
Distribution of the Sum of $n$ Non-Identically Distributed Uniform
Random Variables}

%\date{\today}

%\subjclass{Primary: 60G50; Secondary: 60E10, 42A38}

%\keywords{Uniform distribution, probability density, convolution,
%Fourier transform, sine integrals}

\maketitle

%\begin{center}
%   {\em Running head:} Sum of Non-Identical Uniform Random Variables
%\end{center}

\begin{center}
   {\sc David M. Bradley and Ramesh C. Gupta}

   {\em Department of Mathematics and Statistics, University of
        Maine, Orono, ME 04469-5752, U.S.A.}

\vskip.1in

   {\em e-mail: {\tt dbradley}@{\tt e-math.ams.org},
                {\tt  rcgupta}@{\tt maine.maine.edu} }
\end{center}

\vskip.5in

%\begin{center}
%    (Received: August 28, 2000; accepted: May 16, 2001)
%\end{center}

\noindent{{\em Key words and phrases:} Uniform distribution,
probability density, convolution, Fourier transform, sine
integrals.}

\vskip.1in

%\author{David~M. Bradley}
%\address{Department of Mathematics \& Statistics\\
%         University of Maine\\
%         5752 Neville Hall
%         Orono, Maine 04469-5752\\
%         U.S.A.}
%\email[David~M. Bradley]{bradley@gauss.umemat.maine.edu}
%\author{Ramesh~C. Gupta}
%\email[Ramesh~C. Gupta]{rcgupta@maine.maine.edu}

%\begin{abstract}
{\sc Abstract.}
   The distribution of the sum of independent identically
   distributed uniform random variables is
   well-known.  However, it is sometimes necessary to analyze
   data which have been drawn from different uniform
   distributions.  By inverting the characteristic function, we
   derive explicit formul{\ae} for
   the distribution of the sum of $n$ non-identically
   distributed uniform random variables in both the continuous
   and the discrete case.  The results, though involved, have a
   certain elegance.  As examples, we derive from our general
   formul{\ae} some special cases which have appeared in the literature.
%\end{abstract}

\section{Introduction}\label{sect:Intro}
The classical uniform distribution is perhaps the most versatile
statistical model: applications abound in nonparametric statistics
and Bayes procedures.  Chu (1957) and Leone (1961) utilized
uniform distributions in connection with sample quasi-ranges. Naus
(1966) applied uniform distribution in a power comparison of tests
of random clustering.  For additional applications and examples,
see Johnson et al (1995).

Here, we are concerned with the distribution of the sum of $n$
independent non-identically distributed uniform random variables.
It is well-known that the probability density function of such a
sum, in which the summands are uniformly distributed in a common
interval $[-a,a]$, can be obtained via standard convolution
formul{\ae}:  see Feller (1966, p.\ 27) or Renyi (1970, pp.\
196--197), for example.  However, it is sometimes necessary to
analyze data which have been drawn from non-identical uniform
distributions.  For example, measurements accurate to the nearest
foot may be combined with measurements accurate to the nearest
inch.  In such cases, the distribution of the sum is more
complicated.  Tach (1958) gives tables to five decimal places of
the cumulative distribution of the sum for $n= 2, 3$ and 4 for
some special cases.

The first general result in this direction seems to have been made
by Olds (1952), who derived the distribution of the sum
$\sum_{j=1}^n X_j$, in which each $X_j$ is uniformly distributed
in an interval of the form $[0,a_j)$ with $a_j>0$.  The proof is
by induction, and in that respect is somewhat unsatisfactory,
since in general inductive proofs require knowing beforehand the
formula to be proved.  Subsequently, Roach (1963) deduced what is
essentially Olds' formula using $n$-dimensional geometry.  Later
Mitra (1971), apparently unaware of these previous results,
derived the distribution of the sum in which each random variable
is uniformly distributed in an interval of the form
$[-\omega_j,\omega_j]$ using N\"orlund's (1924) difference
calculus.

Here, we derive an explicit formula for the slightly more general
situation of the distribution of the sum $\sum_{j=1}^n X_j$, in
which each $X_j$ is uniformly distributed in an interval of the
form $[c_j-a_j,c_j+a_j]$ with $a_j>0$ and $c_j$ an arbitrary real
number.  Of course, each of the aforementioned results can be
obtained from ours by specializing the parameters $c_j$ and $a_j$
accordingly.

Our approach is via Fourier theory and is quite straightforward;
specifically we invert the characteristic function.  As a result,
our formula differs somewhat in form from the special cases
alluded to previously.  However, the inversion technique is quite
flexible, and readily lends itself to the study of other types of
distributions, such as the discrete case, which seems not to have
been discussed in the literature. Thus, in a similar fashion, we
derive the distribution of the sum of $n$ random variables with
point mass at the integers in intervals of the form $[-m_j,m_j]$,
in which each $m_j$ is a positive integer. The formula in the
discrete case is somewhat more complicated than the corresponding
formula in the continuous case; nevertheless, they are clearly
closely related, and there is a certain charm and elegance to
both.  Of course, the same results may be obtained using the
standard transformation methods.

\section{The Continuous Case}\label{sect:Cont}
Fix a positive integer $n$, and let $\vec a=(a_1,a_2,\dots,a_n)$
and $\vec c=(c_1,c_2,\dots,c_n)$ be vectors of real numbers with
each $a_j>0$.  For each $j=1,2,\dots,n$, we consider a random
variable $X_j$ uniformly distributed on the closed interval
$[c_j-a_j,c_j+a_j]$.  The step function $\chi_j:\R\to\R$ defined
for real $x$ by
\begin{equation}
   \qquad 2a_j\,\chi_j(x) = \left\{\begin{array}{lll}1 &\mbox{if
   $|x-c_j|<a_j$,}\\ \tfrac12 &\mbox{if $|x-c_j|=a_j$,}\\ 0 &\mbox{if
   $|x-c_j|>a_j$}\end{array}\right.
\label{ChiDef}
\end{equation}
represents the density of the random variable $X_j$ for each
$j=1,2,\dots,n$. (When employing techniques from Fourier theory,
it is convenient to define the densities at jump discontinuities
so that the equation $\chi_j(x+)+\chi_j(x-)=2\chi_j(x)$ is
satisfied for all real $x$.)  The corresponding characteristic
function (Fourier transform) is given by
\[
   \qquad \widehat{\chi}_j(t):=\frac{1}{2a_j}\int_{c_j-a_j}^{c_j+a_j}
   e^{itx}\,dx
   = e^{i c_j t}\,\sinc(a_j t),\qquad t\in\R,
\]
where $\sinc(x):=x^{-1}\sin x$ if $x\ne 0$, and $\sinc(0):=1$. For
real $x$, the density of the sum $\sum_{j=1}^n X_j$ is given by
the $n$-fold convolution $f_n(x) :=
(\chi_1*\chi_2*\dots*\chi_n)(x)$.   Thus,
\[
  \qquad f_n(x) =
  \int_{-\infty}^\infty\chi_1(x-y_2)\int_{-\infty}^\infty
   \chi_2(y_2-y_3)
   \cdots\int_{-\infty}^\infty\chi_{n-1}(y_{n-1}-y_n)
   \chi_n(y_n)\,dy_2\dots dy_n.
\]
In particular, if each of the $n$ intervals is centered at $0$,
and we write $x_1=x-y_2$, $x_n=y_n$ and $x_j=y_j-y_{j+1}$ for
$1<j<n$, then the conditions on the variables $x_1,x_2,\dots,x_n$
are that each $|x_j|<a_j$ and $\sum_{j=1}^n x_j = x.$  Thus,
$f_n(x)$ is simply the volume (in the sense of Lebesgue measure)
of the $(n-1)$-dimensional region
\[
   \qquad\{(x_1,x_2,\dots,x_{n-1})\in\R^{n-1} :
   \big|x-\textstyle\sum_{j=1}^{n-1}x_j\big|<a_n
   \;{\mathrm{and}}\; |x_j|<a_j\;
   {\mathrm{for}}\; 1\le j<n\}
\]
divided by the volume $\prod_{j=1}^n 2a_j$ of the $n$-dimensional
hyperbox
\begin{equation}
   \qquad
   \{(x_1,x_2,\dots,x_n)\in\R^n : |x_j|<a_j \quad {\mathrm{for}}
   \quad j=1,2,\dots,n\}.
\label{hypercube}
\end{equation}
Despite the utility of these representations, it is desirable to
have an explicit formula for $f_n$.  In this vein, we have the
following
\begin{Thm}\label{Thm:Cts}  The density of the sum of
$n$ independent random variables, uniformly distributed in the
intervals $[c_j-a_j,c_j+a_j]$ for $j=1,2,\dots,n$, is given by
\begin{equation}
\label{CtsFormula}
\begin{split}
  \qquad f_n(x) &= \bigg[
  \sum_{\vec\varepsilon\in\{-1,1\}^n}\bigg(x+\sum_{j=1}^n (\varepsilon_j
  a_j-c_j)\bigg)^{n-1}\\
  &\qquad\times\sign\bigg(x+\sum_{j=1}^n (\varepsilon_j
  a_j-c_j)\bigg)\prod_{j=1}^n\varepsilon_j\bigg]
  \bigg/\bigg[(n-1)! 2^{n+1}\prod_{j=1}^n a_j\bigg],
\end{split}
\end{equation}
in which the sum is over all $2^n$ vectors of signs
\[
   \qquad\vec\varepsilon=(\varepsilon_1,\varepsilon_2,\dots,\varepsilon_n)\in
   \{-1,1\}^n \qquad \mbox{i.e. each $\varepsilon_j=\pm 1$}
\]
and
\[
   \qquad
   \sign(y) := \left\{\begin{array}{rrr} 1 &\mbox{if $y>0$,}\\
      0 &\mbox{if $y=0$,}\\
      -1 &\mbox{if $y<0$.}\end{array}\right.
\]
\end{Thm}

\begin{proof} Since the random variables are assumed to be independent,
the characteristic function of the distribution of the sum is the
product of the characteristic functions of their distributions:
\begin{equation}
   \qquad \widehat{f}_n(t) = \prod_{j=1}^n \widehat{\chi}_j(t)
   =\prod_{j=1}^n e^{i tc_j}\,\sinc(a_j t).
\label{sincprod}
\end{equation}
Mitra (1971, p.\ 195) remarks that ``It is difficult to obtain an
inverse Fourier transform of this product in a neat form.''
Indeed, his approach is to expand the product of sincs into a
power series using generalized Bernoulli polynomials. However, we
shall see that the inverse Fourier transform of~(\ref{sincprod})
has an elegant representation involving a sum over the vertices of
the hyperbox~(\ref{hypercube}).

Since for all real $x$, $f_n(x) = \tfrac12 f_n(x+) + \tfrac12
f_n(x-)$ by continuity of $f_n$ for $n>1$, and by definition of
$\chi_1$ when $n=1$, Fourier inversion gives
\[
   \qquad f_n(x) = \frac1{2\pi}\int_{-\infty}^\infty e^{-itx}
   \widehat{f}_n(t)\,dt =
   \frac1{2\pi}\int_{-\infty}^\infty e^{-ity}\prod_{j=1}^n
   \sinc(a_j t)\,dt,
\]
where $y:=x-\sum_{j=1}^n c_j$.  Since $\sinc$ is an even function,
making the change of variable $t\mapsto-t$ yields
\begin{equation}
   \qquad f_n(x) =
 \frac1{2\pi}\int_{-\infty}^\infty e^{ity}\prod_{j=1}^n
   \sinc(a_j t)\,dt.
\label{SincIntegral}
\end{equation}
It remains only to evaluate the integral~(\ref{SincIntegral}).
Related integrals are studied in Borwein and Borwein (2001) using
a version of the Parseval/Plancherel formula, and by expanding the
product of sincs into a sum of cosines.  Our approach is somewhat
more direct.  We first express the $\sinc$ functions using complex
exponentials, so that
\begin{equation}
   \qquad f_n(x) = \frac1{2\pi}\left(\frac1{2i}\right)^n
   \bigg(\prod_{j=1}^n a_j^{-1}\bigg)\int_{-\infty}^\infty t^{-n}
   e^{ity}\prod_{j=1}^n\left(e^{ita_j}-e^{-ita_j}\right)\,dt.
\label{RequiredIntegral}
\end{equation}
For each of the $2^n$ vectors of signs
$\vec\varepsilon=(\varepsilon_1,\varepsilon_2,\dots,\varepsilon_n)\in\{-1,1\}^n$,
let
\[
   \qquad\rho_{\varepsilon}:=\prod_{j=1}^n \varepsilon_j,
   \qquad\mbox{and}\qquad \vec\varepsilon\cdot\vec a=\sum_{j=1}^n
   \varepsilon_j a_j.
\]
By carefully expanding the product of exponentials
in~(\ref{RequiredIntegral}), we find that
\begin{equation}
   \qquad\prod_{j=1}^n \left(e^{ita_j}-e^{-ita_j}\right)
   =\sum_{\vec\varepsilon\in\{-1,1\}^n} \rho_{\varepsilon}
   \exp(it\vec\varepsilon\cdot\vec a).
\label{ProdtoSum}
\end{equation}
It follows that
\begin{equation}
\label{PVIntegral}
   \qquad f_n(x) = \frac1{2\pi}\left(\frac1{2i}\right)^n
   \bigg(\prod_{j=1}^n a_j^{-1}\bigg)^n
   \sum_{\vec\varepsilon\in\{-1,1\}^n}
   \rho_{\varepsilon}\, {\mathrm{P.V.}}\int_{-\infty}^\infty
   t^{-n}\exp(it(y+\vec\varepsilon\cdot\vec a))\,dt.
\end{equation}
Although each of the individual integrals in~(\ref{PVIntegral}) is
divergent, the singularities must cancel
because~(\ref{RequiredIntegral}) converges. Therefore, the
required finite integral~(\ref{RequiredIntegral}) is equal to its
principal value, and hence by linearity is given
by~(\ref{PVIntegral}).

In view of the fact that~(\ref{ProdtoSum}) is entire with a zero
of order $n$ at $t=0$, we may integrate~(\ref{PVIntegral}) by
parts $n-1$ times and thereby obtain
\begin{equation}
\label{PVParts}
\begin{split}
  \qquad  f_n(x) &=
  \frac1{2\pi}\left(\frac1{2i}\right)^n\bigg(\prod_{j=1}^n
  a_j^{-1}\bigg)\frac{i^{n-1}}{(n-1)!}
  \sum_{\vec\varepsilon\in\{-1,1\}^n}\rho_{\varepsilon}
  (y+\vec\varepsilon\cdot\vec a)^{n-1}\\
  & \qquad\times
  {\mathrm{P.V.}}\int_{-\infty}^\infty t^{-1}
  \exp(it(y+\vec\varepsilon\cdot\vec a))\,dt.
\end{split}
\end{equation}
But for any real number $b$, we have
\begin{align*}
  \qquad&{\mathrm{P.V.}}\int_{-\infty}^\infty
  t^{-1}\exp(itb)\,dt\\
  &= \lim_{\varepsilon\to0+}\left\{\int_{\varepsilon}^\infty
  t^{-1}\exp(itb)\,dt
  +\int_{-\infty}^{-\varepsilon}t^{-1}\exp(itb)\,dt\right\}\\
  &= \lim_{\varepsilon\to0+}\left\{\int_{\varepsilon}^\infty
  t^{-1}\exp(itb)\,dt - \int_{\varepsilon}^\infty
  t^{-1}\exp(-itb)\,dt\right\}\\
  &= 2i\int_0^\infty t^{-1}\sin(tb)\,dt\\
  &= i\pi\,\sign(b).
\end{align*}
Applying this latter result to~(\ref{PVParts})
yields~(\ref{CtsFormula}) and completes the proof of
Theorem~\ref{Thm:Cts}.
\end{proof}

In some applications, it may be easier to work with powers of
expressions involving the maximum function $y_+:=\max(y,0)$ as
opposed to the $\sign$ functions in~(\ref{CtsFormula}).  To this
end, we make the following
\begin{Def}\label{TauDef} Let $\tau:\R\to\R$ be given by
\begin{equation}
   \qquad\tau(x) = \left\{\begin{array}{lll} 1,& \mbox{if $x>0$,}\\
   \tfrac12, & \mbox{if $x=0$,}\\ 0, &\mbox{if
   $x<0$,}\end{array}\right.
\label{taudef}
\end{equation}
and for $y$ real and $n$ a positive integer, let $y^{n-1}_+:=
y^{n-1}\tau(y)$.  Note that $y^0_+ = \tau(y)$ and
$y^n_+=(\max(y,0))^n$ for $n>0$.
\end{Def}
Then we have the following corollary to Theorem~\ref{Thm:Cts}.
\begin{Cor}\label{Cor:tau}
Let $f_n$ be as in Theorem~\ref{Thm:Cts}. Then
\begin{equation}
\label{CtsFormula2}
\begin{split}
  \qquad f_n(x) &=
  \bigg[\sum_{\vec\varepsilon\in\{-1,1\}^n}\bigg(x+\sum_{j=1}^n (\varepsilon_j
  a_j-c_j)\bigg)^{n-1}_+\;\prod_{j=1}^n\varepsilon_j\bigg]\bigg/
  \bigg[(n-1)!\,2^n\prod_{j=1}^n a_j\bigg].
\end{split}
\end{equation}
\end{Cor}

\begin{proof}  Note that $\sign(x)=2\tau(x)-1$ holds for all real $x$.
Hence, substituting the $\tau$ function for the $\sign$ function
in~(\ref{CtsFormula}), we see that it suffices to prove the
identity
\begin{equation}
\label{CoolIdentity}
   \qquad\sum_{\vec\varepsilon\in\{-1,1\}^n} \bigg(x+
   \sum_{j=1}^n(\varepsilon_ja_j-c_j)\bigg)^{n-1}\;\prod_{j=1}^n
   \varepsilon_j = 0.
\end{equation}
Since $e^{a_j t}-e^{-a_j t}=2a_j t+O(t^2)$ as $t\to0$,
(\ref{CoolIdentity}) follows easily on comparing coefficients of
$t^{n-1}$ in
\[
   \qquad\sum_{\vec\varepsilon\in\{-1,1\}^n}
   \exp\bigg\{\bigg(x+\sum_{j=1}^n \varepsilon_j a_j-c_j\bigg)t\bigg\}
   \;\prod_{j=1}^n \varepsilon_j
   = e^{xt}\prod_{j=1}^n e^{-c_j
   t}\left(e^{a_jt}-e^{-a_jt}\right).
\]
Alternatively, note that if $x\ge \sum_{j=1}^n (c_j+a_j)$, then
$f_n(x)=0$ by definition:  since each $X_j$ is
$u[c_j-a_j,c_j+a_j]$, the sum $\sum_{j=1}^n X_j$ must fall within
the interval $[\sum_{j=1}^n (c_j-a_j),\sum_{j=1}^n (c_j+a_j)]$.
But, if $x\ge \sum_{j=1}^n (c_j+a_j)$, then we can drop the
subscripted ``+'' from~(\ref{CtsFormula2}) since $x+\sum_{j=1}^n
(\varepsilon_ja_j-c_j)\ge 0$ for each
$\vec\varepsilon\in\{-1,1\}^n$.  It follows
that~(\ref{CoolIdentity}) holds for all $x\ge \sum_{j=1}^n
(c_j+a_j)$.  Since the left hand side of~(\ref{CoolIdentity}) is a
polynomial in $x$ which vanishes for all sufficiently large values
of $x$, it must in fact vanish for all real $x$ by the identity
theorem.
\end{proof}

\begin{Ex}  When $n=1$,
formul{\ae}~(\ref{CtsFormula}) and~(\ref{CtsFormula2}) give the
central difference representations
\[
   \qquad\chi_1(x) = \frac{\sign(x-c_1+a_1)-\sign(x-c_1-a_1)}{4a_1}
    = \frac{\tau(x-c_1+a_1)-\tau(x-c_1-a_1)}{2a_1}
\]
respectively.  Both are equivalent to the
definition~(\ref{ChiDef}) with $j=1$.
\end{Ex}

\begin{Ex}  When $n=2$, we know that
\[
   \qquad f_2(x) =
   (\chi_1*\chi_2)(x)=\int_{-\infty}^\infty
   \chi_1(x-y)\chi_2(y)\,dy.
\]
If $c_1=c_2=0$, the convolution reduces to
\[
  \qquad f_2(x)
   =\frac{1}{4a_1a_2}\int_{\max(x-a_1,\,-a_2)}^{\min(x+a_1,\,a_2)}\,dy
   =\frac{\min(x+a_1,a_2)-\max(x-a_1,-a_2)}{4a_1a_2},
\]
so that, in particular,
\begin{equation}
   \qquad f_2(0)
   = \frac{1}{2\pi}\int_{-\infty}^\infty
   \frac{\sin(a_1t)}{a_1t}\cdot\frac{\sin(a_2t)}{a_2t}\,dt =
   \frac{\min(a_1,a_2)}{2a_1a_2}. \label{Checkn=2}
\end{equation}
On the other hand, in light of the fact that $y\,\sign(y)=|y|$ for
$y$ real, formula~(\ref{CtsFormula}) gives
\begin{align*}
   \qquad f_2(0) &=
   \frac{|a_1+a_2|-|a_1-a_2|-|-a_1+a_2|+|-a_1-a_2|}{8a_1a_2}\\
    &= \frac{2(a_1+a_2)-2|a_1-a_2|}{8a_1a_2}\\
    &= \frac{\min(a_1,a_2)}{2a_1a_2},
\end{align*}
in agreement with~(\ref{Checkn=2}).  Alternatively, since
$y\tau(y)=y_+=\max(y,0)$, formula~(\ref{CtsFormula}) gives
\begin{align*}
   \qquad f_2(0) &=
   \frac{\max(a_1+a_2,0) -\max(a_1-a_2,0)-\max(-a_1+a_2,0)+\max(-a_1-a_2,0)}
   {4a_1a_2}\\
   &= \frac{ a_1+a_2-\max(a_1-a_2,0)-\max(a_2-a_1,0)}{4a_1a_2}\\
   &= \frac{\min(a_1,a_2)}{2a_1a_2}.
\end{align*}
\end{Ex}

\begin{Ex} If the random variables comprising the sum
are uniformly distributed in a common interval centered at 0---say
$[-a,a]$ with $a>0$---then formula~(\ref{CtsFormula2}) gives
\[
   \qquad f_n(x) =
   \frac1{(n-1)!\,(2a)^n}
   \sum_{\vec\varepsilon\in\{-1,1\}^n}
   \bigg(x+a\sum_{j=1}^n\varepsilon_j\bigg)^{n-1}_+\;
   \prod_{j=1}^n\varepsilon_j.
\]
If $k$ of $\varepsilon_1,\varepsilon_2,\dots,\varepsilon_n$ are
negative and the remaining $n-k$ are positive, then summing over
$k$ yields
\[
   \qquad f_n(x) = \frac1{(n-1)!\,(2a)^n}\sum_{k=0}^n (-1)^k \binom{n}{k}
      \left(x+(n-2k)a\right)^{n-1}_+,
\]
in agreement with Feller (1966).

\end{Ex}

\begin{Ex}  If for
each $j=1,2,\dots,n$ we set $c_j=a_j$ and then replace $a_j$ by
$a_j/2$, then $X_j$ will be uniformly distributed in $[0,a_j]$ and
will have density $\chi_j$ now given by
\begin{equation}
\label{NewChiDef}
   \qquad a_j\,\chi_j(x) :=\left\{\begin{array}{lll}
   1, &\mbox{if $0<x<a_j$,}\\
   \tfrac12, &\mbox{if $x=a_j$ or if $x=0$,}\\
   0, &\mbox{if $x>a_j$ or if $x<0$.}\end{array}\right.
\end{equation}
Formula~(\ref{CtsFormula2}) of Corollary~\ref{Cor:tau} now gives
\[
\begin{split}
   \qquad f_n(x) &=\bigg[
   \sum_{\vec\varepsilon\in\{-1,1\}^n}\bigg(x-\sum_{j=1}^n \bigg(
   \frac{1-\varepsilon_j}{2}\bigg)a_j\bigg)^{n-1}_+\;
   \prod_{j=1}^n\varepsilon_j\bigg]
   \bigg/\bigg[(n-1)!\,\prod_{j=1}^n a_j\bigg]\\
   &=\bigg[\sum_{\vec s\in\{0,1\}^n} (-1)^{\Sigma \vec s}
   \left(x-\vec s\cdot\vec a\right)^{n-1}_+\bigg]
   \bigg/\bigg[(n-1)!\,\prod_{j=1}^n a_j\bigg],
\end{split}
\]
where the sum is now over all $2^n$ vectors $\vec
s=(s_1,s_2,\dots,s_n)$ in which each component takes the value $0$
or $1$, $\Sigma \vec s$ denotes the sum of the components
$s_1+s_2+\cdots +s_n$ and $\vec s\cdot \vec a$  denotes the dot
product $s_1 a_1+s_2a_2+\cdots+s_na_n$.  If we now break up the
sum according to the number of non-zero components in the vector
$\vec s$, we find that
\[
\begin{split}
\label{Olds}
   \qquad f_n(x) &=
   \bigg[ x^{n-1}_+ - \sum_{1\le j_1\le n}
   \left(x-a_{j_1}\right)^{n-1}_+ +\sum_{1\le j_1 < j_2\le n}
   \left(x-a_{j_1}-a_{j_2}\right)^{n-1}_+ -+\cdots\\
   &\qquad\qquad+(-1)^n
   \sum_{1\le j_1 < j_2 <\cdots < j_n\le n}
   \bigg(x-\sum_{k=1}^n a_{j_k}\bigg)^{n-1}_+\,\bigg]
    \bigg/\bigg[(n-1)!\,\prod_{j=1}^n a_j\bigg],
\end{split}
\]
which is equivalent to the formula of Olds (1952, p.\ 282 (1))
when $n>1$.  When $n=1$, the two formul{\ae} differ at $x=0$ and
at $x=a_1$ because Olds uses the density which is $1/a_1$ for
$0\le x<a_1$ and $0$ otherwise, in contrast to our more
symmetrical density~(\ref{NewChiDef}).
\end{Ex}

\section{The Discrete Case}\label{sect:Discrete}
Here, we fix $n$ positive integers $m_1,m_2,\dots,m_n$.  For each
$j=1,2,\dots,n$, we now consider a random variable $X_j$ uniformly
distributed on the set of $(2m_j+1)$ integers contained in the
closed interval $[-m_j,m_j]$, i.e.\ the set
$\{-m_j,1-m_j,\dots,m_j-1,m_j\}$. The mass function of $X_j$ is
the rational-valued function of an \emph{integer} variable given
by
\begin{equation}
   \qquad(2m_j+1)\chi_j(p) := \left\{\begin{array}{ll}1 &\mbox{if $|p|\le
   m_j$,}\\ 0 &\mbox{if $|p|>m_j$.}\end{array}\right.
\label{newChiDef}
\end{equation}
We seek a formula analogous to~(\ref{CtsFormula}) for the
probability mass function of the sum $\sum_{j=1}^n X_j$, namely
the $n$-fold convolution
\begin{equation}
   \qquad g_n(p) := (\chi_1*\chi_2*\cdots*\chi_n)(p)
   = \sum_{k_1+k_2+\cdots+k_n=p}\; \prod_{j=1}^n \chi_j(k_j),
\label{DiscreteConvolution}
\end{equation}
where the sum is over all integers $k_1,k_2,\dots,k_n$ such that
$k_1+k_2+\cdots+k_n=p$.  Although~(\ref{DiscreteConvolution}) is a
formula of sorts, the various conditions on the summation indices
make it inconvenient to apply.  For example, when $n=3$, let
$u_k=u_k(p):=\min(m_2,p-k+m_3)$ and
$v_k=v_k(p):=\max(-m_2,p-k-m_3)$. Then
\[
    \qquad g_3(p) =\bigg(\prod_{j=1}^3 (2m_j+1)^{-1}\bigg)
    \sum_{\substack{|k|\le m_1\\u_k<v_k}}
    (u_k-v_k+1),
\]
where the sum is over all integers $k$ for which $u_k<v_k$ and
$|k|\le m_1$. In general, one needs to consider cases which depend
on the size of $p$ in relation to various signed sums of subsets
of the parameters $m_1,m_2,\dots,m_n$.  The number of cases to be
delineated increases exponentially with $n$. Thus, the situation
becomes rapidly unwieldy. Fortunately, there is alternative
approach provided by Fourier theory.  With the convolution
$g_n:\Z\to\Z$ defined as above, we have
\begin{Thm}\label{Thm:Discrete} For all integers $p$,
\begin{equation}
\begin{split}
\label{DiscreteFormula}
   \qquad g_n(p) &= \frac{M}{2^n}\sum_{\vec\varepsilon\in\{-1,1\}^n}
   \sign\bigg(2p+\sum_{j=1}^n (2m_j+1)\varepsilon_j\bigg)
   \bigg(\prod_{j=1}^n\varepsilon_j\bigg)\\
   &\qquad\qquad\times\sum_{k=0}^{(n-1)/2}
   (-1)^k
   b_{2k}^{(n)}
   \frac{\left(2p+\sum_{j=1}^n(2m_j+1)\varepsilon_j\right)^{n-2k-1}}
   {(n-2k-1)!},
\end{split}
\end{equation}
where
\[
   \qquad M:=\prod_{j=1}^n (2m_j+1)^{-1},
\]
and the rational numbers $b_{2k}^{(n)}$ are the coefficients in
the Laurent series expansion
\begin{equation}
   \qquad \bigg(\frac1{\sin x}\bigg)^n = \sum_{k=0}^\infty x^{2k-n}
    b_{2k}^{(n)}.
\label{CscLaurent}
\end{equation}
Explicitly,
\begin{equation}
   \qquad b_{2k}^{(n)} = (-1)^k \binom{n+2k}{n}\sum_{m=0}^{2k}
   \frac{n}{n+m}\binom{2k}{m}\frac{1}{2^m(2k+m)!}
   \sum_{r=0}^m (-1)^r\binom{m}{r}(2r-m)^{2k+m}.
\label{LaurentCoefficientFormula}
\end{equation}
\end{Thm}

\begin{proof}
The corresponding characteristic functions for $j=1,2,\dots,n$ are
now given by
\[
   \qquad(2m_j+1)\widehat{\chi}_j(t) = \sum_{|k|\le m_j} e^{itk}
   = \left\{\begin{array}{ll}
     \displaystyle{\frac{\sin((m_j+1/2)t)}{\sin(t/2)}},
     &\mbox{if $t\notin 2\pi\Z$,}\\ 2m_j+1, &\mbox{if $t\in
     2\pi\Z$,}\end{array}\right.
\]
which we recognize as the familiar Dirichlet kernel---see
eg.~Korner (1988, p.\ 68). Since for $t\notin 2\pi\Z$, we have
\[
   \qquad\left(\chi_1*\chi_2*\cdots*\chi_n\right)\sphat (t)
   =\prod_{j=1}^n \widehat{\chi}_j(t)
   =M\prod_{j=1}^n \frac{\sin((m_j+1/2)t)}{\sin(t/2)},
\]
it follows by orthogonality of the exponential that
\[
   \qquad g_n(p) = \frac{M}{2\pi}\int_0^{2\pi} e^{ipt}\prod_{j=1}^n
   \frac{\sin((m_j+1/2)t)}{\sin(t/2)}\,dt
   = \frac{M}{\pi}\int_0^\pi e^{2ipx}\prod_{j=1}^n
   \frac{\sin((2m_j+1)x)}{\sin x}\,dx.
\]
We remark in passing that the factors $\sin((2m_j+1)x)/\sin x$ in
this latter representation are simply the Chebyshev polynomials
$U_{2m_j}(\cos x)$ of the second kind---see eg.~Abramowitz and
Stegun (1972, p.\ 766). The partial-fraction expansion
\begin{equation}
   \qquad\bigg(\frac{1}{\sin x}\bigg)^n
   = \sum_{k=0}^{(n-1)/2} b_{2k}^{(n)}
   \sum_{r=-\infty}^\infty\bigg(\frac{(-1)^r}{x+r\pi}\bigg)^{n-2k},
\label{CscParfrac}
\end{equation}
is a modified version of the formula in Schwatt (1924, pp.\
209--210), and yields
\begin{equation}
\label{GnParfrac}
   \qquad g_n(p) = \sum_{k=0}^{(n-1)/2} b_{2k}^{(n)}
   \sum_{r=-\infty}^\infty\frac{M}{\pi}
   \int_0^\pi e^{2ipx}\bigg(\frac{(-1)^r}{x+r\pi}\bigg)^{n-2k}
   \prod_{j=1}^n \sin((2m_j+1)x)\,dx.
\end{equation}
When $n-2k>1$, the bilateral series
\[
   \qquad \sum_{r=-\infty}^\infty \bigg(\frac{(-1)^r}{x+r\pi}\bigg)^{n-2k}
\]
converges absolutely, and the interchange of summation and
integration is easily justified using either Lebesgue's dominated
convergence theorem or Fubini's theorem with both Lebesgue and
counting measure. If $n-2k=1$,  absolute convergence can be
recovered by recasting the bilateral series in the form
\[
   \qquad \sum_{r=-\infty}^\infty \frac{(-1)^r}{x+r\pi}
   = \frac1x+\sum_{r=1}^\infty (-1)^r\bigg(\frac{ 2 x}{x^2-r^2\pi}\bigg),
\]
and the justification proceeds as in the previous case.

We shall see that the inner sum of integrals in~(\ref{GnParfrac})
can be expressed as a single integral over the whole real line. To
this end, we compute
\begin{align*}
   \qquad &\sum_{r=-\infty}^\infty\frac1{\pi}
   \int_0^\pi e^{2ipx}\bigg(\frac{(-1)^r}{x+r\pi}\bigg)^{n-2k}
   \prod_{j=1}^n \sin((2m_j+1)x)\,dx\\
  =&\sum_{r=-\infty}^\infty\frac{(-1)^{(n-2k)r}}{\pi}
   \int_{r\pi}^{(r+1)\pi}
   t^{2k-n} e^{2ip(t-r\pi)}\prod_{j=1}^n
   \sin\left((2m_j+1)(t-r\pi)\right)\,dt\\
   =&\sum_{r=-\infty}^\infty\frac{(-1)^{(n-2k)r}}{\pi}
   \int_{r\pi}^{(r+1)\pi} t^{2k-n}e^{2ipt}\prod_{j=1}^n (-1)^{(2m_j+1)r}
   \sin((2m_j+1)t)\,dt\\
   =&\;\;\frac1{\pi}\int_{-\infty}^\infty t^{2k-n} e^{2ipt}\prod_{j=1}^n
   \sin((2m_j+1)t)\,dt.
\end{align*}
This latter integral can be evaluated just as we evaluated the
integral~(\ref{SincIntegral}).  After expanding the product of
sines as a sum over the constituent exponentials and integrating
by parts $n-2k-1$ times, one finds that
\begin{multline}
   \frac1{\pi}\int_{-\infty}^\infty t^{2k-n} e^{2ipt}\prod_{j=1}^n
   \sin((2m_j+1)t)\,dt\\
    =
   \frac{(-1)^k}{2^n(n-2k-1)!}
   \sum_{\vec\varepsilon\in\{-1,1\}^n}
   \bigg(2p+\sum_{j=1}^n(2m_j+1)\varepsilon_j
   \bigg)^{n-2k-1}\\
   \times\sign\bigg(2p+\sum_{j=1}^n(2m_j+1)\varepsilon_j\bigg)
   \prod_{j=1}^n\varepsilon_j.
\label{latterformula}
\end{multline}
Substituting~(\ref{latterformula}) into~(\ref{GnParfrac}) and
interchanging the order of summation completes the proof of
Theorem~\ref{Thm:Discrete}.
\end{proof}

\noindent{\bf Remark.} The equations~(\ref{CscLaurent}),
(\ref{LaurentCoefficientFormula}), (\ref{CscParfrac}) can also be
obtained using Jordan's (1979, \S74, p.\ 216) general formula~ for
the higher derivatives of a power of a reciprocal function and
then applying Mittag-Leffler's theorem.

As in the continuous case, we can replace the $\sign$ function in
Theorem~\ref{Thm:Discrete} with the $\tau$ function of
Definition~\ref{TauDef}.  Thus we obtain
\begin{Cor}\label{Cor:DiscreteFormula2} Let $g_n$, $M$, and
$b_{2k}^{(n)}$ be as in Theorem~\ref{Thm:Discrete}.  Then, for all
positive integers $n$ and integer $p$,
\begin{equation}
\label{DiscreteFormula2}
    \qquad g_n(p) = \frac{M}{2^{n-1}}
    \sum_{k=0}^{(n-1)/2} \frac{(-1)^k b_{2k}^{(n)}}{(n-2k-1)!}
    \sum_{\vec\varepsilon\in\{-1,1\}^n}
    \bigg(2p+\sum_{j=1}^n(2m_j+1)\varepsilon_j\bigg)^{n-2k-1}_+\;
    \prod_{j=1}^n\varepsilon_j,
\end{equation}
where $y^{n-2k-1}_+=y^{n-2k-1}\tau(y)$ as in
Definition~\ref{TauDef}.
\end{Cor}

\begin{proof}  Making the substitution $\sign(y)=2\tau(y)-1$ in
Theorem~\ref{Thm:Discrete} and noting that the coefficient of
$t^{n-2k-1}$ in
\[
   \qquad\sum_{\vec\varepsilon\in\{-1,1\}^n}\exp\bigg\{\bigg(
   2p+\sum_{j=1}^n(2m_j+1)\varepsilon_j\bigg)t\bigg\}=
   e^{2pt}\prod_{j=1}^n\left(e^{(2m_j+1)t}-e^{-(2m_j+1)t}\right)
\]
vanishes for $0\le k\le (n-1)/2$, we obtain the desired result.
\end{proof}

\begin{Ex}  When $n=1$, formul{\ae}~(\ref{DiscreteFormula})
and~(\ref{DiscreteFormula2}) give the representations
\[
\begin{split}
   \qquad(2m_1+1)\chi_1(p) &= \tfrac12\sign(p+m_1+\tfrac12)
                    -\tfrac12\sign(p-m_1-\tfrac12)\\
   &= \tau(p+m_1+\tfrac12)-\tau(p-m_1-\tfrac12),
\end{split}
\]
respectively.  Both are equivalent to the
definition~(\ref{newChiDef}) with $j=1$.
\end{Ex}

\begin{Ex}  When $n=2$, we have $M:=(2m_1+1)^{-1}(2m_2+1)^{-1}$, and
\[
   \qquad g_2(p) =
   \chi_1*\chi_2(p)= \sum_{k+j=p}\chi_1(k)\chi_2(j).
\]
Since this is simply the coefficient of $t^p$ in the product
\[
   \qquad M   \sum_{|k|\le m_1} t^k \sum_{|j|\le m_2} t^j,
\]
letting $u(p):=\min(m_1,p+m_2)$ and $v(p):=\max(-m_1,p-m_2)$, we
have
\[
   \qquad g_2(p) = M\times\left\{\begin{array}{ll} u(p)-v(p)+1 &\mbox{if
   $u(p)\ge v(p)$,}\\ 0 &\mbox{if $u(p)<v(p)$.}\end{array}\right.
\]
On the other hand formula~(\ref{DiscreteFormula}) gives the
elegant representation
\begin{align*}
   \qquad g_2(p) &=
   \tfrac12 M\left(|p+m_1+m_2+1|
   -|p+m_1-m_2|
   -|p-m_1+m_2|\right.\\
   &\left.\qquad\qquad+|p-m_1-m_2-1|\right),
\end{align*}
in which we have used the relation $y\,\sign(y)=|y|$ for $y$ real.
\end{Ex}

%\begin{thebibliography}{99}

\begin{center} {\sc References}
\end{center}

%\bibitem{AS}
\noindent Abramowitz, M.\ and Stegun, I.\ A.\ (1972). {\it
Handbook of Mathematical Functions}, Dover Publications, New York.
\vskip.2in

%\bibitem{BB}
\noindent Borwein, D.\ and Borwein, J.\ M.\ (2001). Some
remarkable properties of sinc and related integrals. {\it The
Ramanujan Journal}, \textbf{5}, (1), 73--89.  \vskip.2in

%\bibitem{Chu57}
\noindent Chu, J.\ T.\ (1957). Some uses of quasi-ranges, {\it
Annals of Mathematical Statistics}, \textbf{28}, 173--180.
\vskip.2in

%\bibitem{Fell}
\noindent Feller, W.\ (1966). {\it An Introduction to Probability
Theory and its Applications}, Vol.\ II, John Wiley \& Sons, New
York.\vskip.2in

%\bibitem{John95}
\noindent Johnson, N.\ L., Kotz, S.\ and Balakrishnan, N.\ (1995).
{\it Continuous univariate distributions}, (2nd ed.) John Wiley,
New York.\vskip.2in

%\bibitem{Jordan}
\noindent Jordan, C.\ (1979). {\it Calculus of Finite
Differences}, Chelsea Publishing, (3rd ed.\ reprinted). \vskip.2in

%\bibitem{Korner}
\noindent K\"orner, T.\ W.\ (1988). {\it Fourier Analysis},
Cambridge University Press.\vskip.2in

%\bibitem{Leon61}
\noindent Leon, F.\ C.\ (1961). The use of sample quasi-ranges in
setting confidence intervals for the population standard
deviation, \textit{Journal of the American Statistical
Association}, \textbf{56}, 260--272.\vskip.2in

%\bibitem{Mitra}
\noindent Mitra, S.\ K.\ (1971). On the probability distribution
of the sum of uniformly distributed random variables, {\it SIAM
J.\ Appl.\ Math.}, \textbf{20}, (2), 195--198.\vskip.2in

%\bibitem{Naus66}
\noindent Naus, I.\ (1966).  A power comparison of two tests of
non-random clustering, {\it Technometrics}, \textbf{8}, 493--517.
\vskip.2in

%\bibitem{Nor}
\noindent  N\"orlund, N.\ E.\ (1924). {\it Vorlesungen \"uber
Differenzenrechnung}, Springer, Berlin.\vskip.2in

%\bibitem{Olds52}
\noindent Olds, E.\ G.\ (1952).  A note on the convolution of
uniform distributions, {\it Annals of Mathematical Statistics},
\textbf{23}, 282--285.\vskip.2in

%\bibitem{Renyi}
\noindent  R\'enyi, A.\ (1970). {\it Probability Theory},
North-Holland Publishing, Amsterdam.\vskip.2in

%\bibitem{Roach63}
\noindent Roach, S.\ A.\ (1963). The frequency distribution of the
sample mean where each member of the sample is drawn from a
different rectangular distribution, {\it Biometrika}, \textbf{50},
508--513.\vskip.2in

%\bibitem{Schwatt}
\noindent  Schwatt, I.\ J.\ (1924). {\it An Introduction to the
Operations with Series}, University of Pennsylvania Press,
Philadelphia.\vskip.2in

%\bibitem{Tach58}
\noindent Tach, L.\ T.\ (1958).  Tables for cumulative
distribution function of a sum of independent random variables,
{\it Convair Aeronautics Report M}, ZU-7-119-TN, San Diego.

%\end{thebibliography}
\end{document}